\documentclass[11pt,reqno]{article}
\usepackage{graphicx,epsfig}
\usepackage{amsmath}

\title {Coexistence for Richardson type competing spatial growth models}

\author{Christopher Hoffman}

\usepackage{latexsym}
\usepackage{amsfonts}
\usepackage{amsmath}
\usepackage{amsthm}
\usepackage{graphicx}

\newtheorem{thm}{Theorem}

\newtheorem{lemma}{Lemma}
\newtheorem{cor}{Corollary}

\newcommand{\be}{\begin{equation}}
\newcommand{\ee}{\end{equation}}

\newcommand{\E}{\mbox{$\bf E$}}
\newcommand{\booze}{Busemann }

\newcommand{\N}{{\mathbb N}}

\newcommand{\zd}{{\mathbb Z}^d}

\newcommand{\R}{{\mathbb R}}

\newcommand{\0} {{\bf 0}}
\newcommand{\ez}{{\bf 1}}

\newcommand{\olle}{{H\"{a}ggstr\"{o}m }}
\newcommand{\fpp}{{\mu}}
\newcommand{\mug}{\mbox{C}}
\newcommand{\edges}{\mbox{Edges}}

\begin{document}

\maketitle

\begin{abstract}
We study a large family of competing spatial growth models.  In
these the vertices in $\zd$ can take on three possible states
\{0,1,2\}.  Vertices in states 1 and 2 remain in their states
forever, while vertices in state 0 which are adjacent to a vertex
in state 1 (or state 2) can switch to state 1 (or state 2). We
think of the vertices in states 1 and 2 as infected with one of
two infections while the vertices in state 0 are considered
uninfected. In this way these models are variants of the
Richardson model. We start the models with a single vertex in
state 1 and a single vertex is in state 2. We show that with
positive probability state 1 reaches an infinite number of
vertices and state 2 also reaches an infinite number of vertices.
This extends results and proves a conjecture of \olle and Pemantle
\cite{HP}. The key tool is applying the ergodic theorem to
stationary first passage percolation.
\end{abstract}

\section{First Passage Percolation}
In this paper we study a class of competing spatial growth models
by first studying stationary first passage percolation and then
applying our results to the spatial growth models. In first
passage percolation every edge in a graph is assigned a
non-negative number.  This is interpreted as the time it takes to
move across the edge.  This model was introduced by Hammersley and
Welsh \cite{HW}.  See \cite{K} for an overview of first passage
percolation.

Let $\fpp$ be a stationary measure on $[0,\infty)^{\edges(\zd)}$
and let $\omega$ be a realization of $\fpp$. For any $x$ and $y$
we define the {\bf passage time from $x$ to $y$}, $\tau(x,y)$, by
$$\tau(x,y)=\inf\sum{\omega(v_{i},v_{i+1})}$$
where the sum is taken over all of the edges in the path and the
$\inf$ is taken over all paths connecting $x$ to $y$.

The most basic result from first passage percolation is the shape
theorem. We let $\0=(0, \dots,0)$ and $\ez=(1,0, \dots,0).$ Define
$$S(t)=\{x:\ \tau(\0,x)\leq t\}$$ and
$$\bar S(t)=S(t)+\left[-\frac12,\frac12\right]^d.$$  The shape theorem says that there is a
 nonempty set $S$ such that $\frac{\bar S(t)}{t}$ converges to
$S$ a.s.
\begin{thm} \cite{B}
\label{shape1}
 Let $\fpp$ be stationary and ergodic, with the distribution on any edge have finite $d+\epsilon$
 moment for $\epsilon>0$.
 There exists a set $S$ which is nonempty, convex, and symmetric
 about reflection through the origin
 such that for every $\epsilon>0$ there exists a $T$ such that
 for all $t>T$
$$P\left((1-\epsilon)S<\frac{\bar S(t)}{t}<(1+\epsilon)S\right)>1-\epsilon.$$
\end{thm}
This theorem is a consequence of Kingman's subadditive ergodic
theorem.  It is the only property of first passage percolation
that we need.  In general little is known about the shape of $S$
other than it is convex and symmetric. Cox and Durrett have shown
that there are nontrivial product measures such that the boundary
of $S$ contains a flat piece \cite{DL}. However for any compact
nonempty convex set $S$ there exist a stationary measure $\fpp$
such that the shape for $\fpp$ is $S$ \cite{H}.

Another widely studied aspect of first passage percolation are
geodesics. A {\bf geodesic} is a path $G=\{v_0,v_1,\dots\}$ such
that
 $$\tau(v_m,v_n)=\sum_{i=m}^{n-1}{\omega(v_{i},v_{i+1})}$$ for any
$m<n$. We let $G^{\omega}(x,y)=G(x,y)$ be the union of all
geodesics connecting $x$ and $y$. Define $$\Gamma(x)=\cup_{y \in
\zd}
    \{e:e\in \edges(\zd)\mbox{ and }e \in G(x,y)\}.$$
We refer to this as the {\bf tree of infection of $x$}. We define
$K(\Gamma(x))$ to be the number of topological ends in
$\Gamma(x)$.  This is also the number of infinite self avoiding
paths in $\Gamma(x)$ that start at $x$.



Newman has conjectured that for a large class of $\fpp$,
$|K(\Gamma(\0))|=\infty$ a.s. \cite{N}  \olle and Pemantle proved
that if $d=2$ and $\fpp$ is the i.i.d\ with exponential
distribution then with positive probability $|K(\Gamma(\0))|>1$.
Newman has proved that if $\mu$ is i.i.d.\ and $S$ has  certain
properties then $|K(\Gamma(\0))|=\infty$ a.s. \cite{N} Although
these conditions are plausible there are no known measures $\fpp$
with $S$ that satisfy these conditions.

Now we will introduce some more notation which will let us list
the conditions that we place on $\fpp$.  We say that the
configuration $\omega$ has {\bf unique geodesics} if for all $x,y
\in \zd$ there exists a unique geodesic from $x$ to $y$.
If there exists a unique geodesic between $x$ and $y$ we denote it
by $G(x,y)$. The configuration $\omega$ has {\bf unique passage
times} for all $x$ and $y\neq z$
$$\tau(x,y)\neq \tau(x,z).$$
For any $\omega$ we let $\fpp^{(\0,\ez)}_{\omega}$ be the
conditional distribution of $\fpp$ on the edge $(\0,\ez)$ given
that $\omega'(v,w)=\omega(v,w)$ for all edges except $(\0,\ez)$.
We say that $\fpp$ has {\bf finite energy} if for any set
$A\subset \R$ such that $\fpp\{\omega(\0,\ez) \in A\}>0$ and
almost every $\omega$,
 $\fpp^{(\0,\ez)}_{\omega}\{\omega':\omega'(0,1)\in A\}>0$.

As $\fpp$ is a stationary measure we can study its ergodic
theoretical properties. For any $v \in \zd$ define the shift map
$T^v:[0,\infty)^{\edges(\zd)} \to [0,\infty)^{\edges(\zd)}$ by
$$T^v(\omega)(j)=\omega(j+v)$$
for all $j \in \mbox{Edges}(\zd)$. The measure $\fpp$ is {\bf
totally ergodic} if for all $v \in \zd$ the action $(\fpp,T^v)$ is
ergodic.

Now we are ready to define the class of measures that we will work
with. We say that $\fpp$ is {\bf good} if
\begin{enumerate}
\item $\fpp$ is totally ergodic,
\item $\fpp$ has all the symmetries of $\zd$, \label{cont0}
\item the distribution of $\fpp$ on any edge has finite $d+\epsilon$ moment for some
      $\epsilon>0$
\item $\fpp$ has finite energy \label{cont1}
\item $\fpp^{(\0,\ez)}_{\omega}$ is an absolutely continuous measure with support $[0,\infty)$ a.s., and \label{cont2}
\item $\fpp$ produces a shape $S$ which is bounded.
\end{enumerate}
Note that conditions \ref{cont0}, \ref{cont1} and \ref{cont2}
imply that $\fpp$ has unique geodesics and unique passage times.
These conditions were chosen to make the arguments as easy as
possible and could be made more general. All that is essential for
the argument to show that there are at least two disjoint infinite
geodesics is that $\fpp$ is totally ergodic and that Lemma
\ref{lem0} and Corollary \ref{shape2} below are satisfied. The
conditions \ref{cont0}, \ref{cont1} and \ref{cont2} are used to
show that coexistence occurs with positive probability. Throughout
the rest of the paper we will assume that $\mu$ is good.
Unfortunately there is no general necessary and sufficient
condition to determine when the shape $S$ is bounded. See \cite{H}
for examples.

\section{Spatial Growth Models}
Now we explain the relationship between first passage percolation
and our competing growth models. For any $\omega \in
[0,\infty)^{\edges(\zd)}$ with unique passage times and any $x\neq
y \in \zd$ we can project it to $\tilde{\omega}_{x,y} \in
\left(\{0,1,2\}^{\zd}\right)^{[0,\infty)}$ by
$$\tilde{\omega}_{x,y}(z,t)=\left\{%
\begin{array}{ll}
     2& \hbox{if $\tau(x,z)\leq t$ and $\tau(x,z) <\tau (y,z)$;} \\
     1& \hbox{if $\tau(y,z)\leq t$ and $\tau(x,z) >\tau (y,z)$;} \\
     0& \hbox{else.} \\
\end{array}%
\right.
$$
If $\fpp$ has unique passage times a.s.\ then $\fpp$ projects onto
a measure on $\left(\{0,1,2\}^{\zd}\right)^{[0,\infty)}$. It is
clear that the models start with a single vertex in state 1 and a
single vertex is in state 2.  Vertices in states 1 and 2 remain in
their states forever, while vertices in state 0 which are adjacent
to a vertex in state 1 (or state 2) can switch to state 1 (or
state 2). We think of the vertices in states 1 and 2 as infected
with one of two infections while the vertices in state 0 are
considered uninfected. In this way these models are variants of
the Richardson model.

As each $z \in \zd$ eventually changes to state 1 or 2 and then
stays in that state for the rest of time, we can speak of the
limiting configuration. There are two possible outcomes. The first
is coexistence or mutual unbounded growth. If this occurs then the
limiting configuration has infinitely many $z$ in state 1 and
infinitely many $z$ in state 2.  The other outcome is domination.
If this happens then in the limiting configuration there are only
finitely many vertices in that state and all but finitely many
vertices are in the other state.

For many measures $\fpp$ (for example if $\fpp$ is i.i.d. with
nontrivial marginals) then it is easy to prove that domination
occurs with positive probability.  But it is much more difficult
to show that coexistence occurs with positive probability. More
precisely we define $\mug(x,y)$ to be the event that
$$|\{z:\lim_{t \to \infty}\tilde{\omega}_{x,y}(z)=1\}|
    =|\{z:\lim_{t \to \infty}\tilde{\omega}_{x,y}(z)=2\}|
    =\infty.$$
We refer to this event as {\bf coexistence} or {\bf mutual
unbounded growth}. Our main result is that with positive
probability coexistence occurs.
\begin{thm} \label{main} If $\mu$ is good then
 $$P(\mug(\0,\ez))>0.$$
\end{thm}
This proves a conjecture of \olle and Pemantle  \cite{HP}. They
proved this theorem in the case that $d=2$ and $\fpp$ is i.i.d.\
with exponential distribution.  Garet and Marchand have given a
different proof of Theorem \ref{main} \cite{G}.  Their method
follows more closely the approach taken by \olle and Pemantle.

\section{Outline}
In this section we outline the proof of our main result. For any
$x,y \in \zd$ and infinite geodesic $G=(v_0,v_1,v_2,\dots)$  we
can define
$$B^{\omega}_{G}(x,y)=B_G(x,y)=\lim_{n \to \infty} \tau(x,v_n)-\tau(y,v_n).$$
To see the limit exists first note that
\begin{eqnarray*}
B_G(x,y)& =& \lim_{n \to \infty} \tau(x,v_n)-\tau(y,v_n)\\
        & =& \lim_{n \to \infty} \tau(x,v_n)-\tau(v_0,v_n)+\tau(v_0,v_n)-\tau(y,v_n)\\
        & =& \lim_{n \to \infty} (\tau(x,v_n)-\tau(v_0,v_n))+\lim_{n \to \infty}(\tau(v_0,v_n)-\tau(y,v_n)).
\end{eqnarray*}
As $G$ is a geodesic the two sequences in the right hand side of
the last line are bounded and monotonic so they converge. Thus
$B_G(x,y)$ is well defined.  If for a given $\omega$ and all $x,y
\in \zd$ the function $B_G(x,y)$ is independent of the choice of
infinite geodesic $G$ then we can define the \booze function
$$B(x,y)=B^{\omega}(x,y)=B^{\omega}_G(x,y).$$
The main step in our proof is Lemma \ref{coalesce}, which states
that the probability that $\{B(x,y)\}_{x,y \in \zd}$ is well
defined is 0.


We will work by contradiction to prove Lemma \ref{coalesce}. In
Lemmas \ref{lem4} and \ref{lem1} we assume that $\{B(x,y)\}_{x,y
\in \zd}$ is well defined a.s.\ and then apply the ergodic theorem
to $\{B(x,y)\}_{x,y \in \zd}$. Then in Lemma \ref{coalesce} we
show that the conclusions of Lemma \ref{lem1} generate a
contradiction with the shape theorem. Thus with positive
probability there are vertices $x$ and $y$ and distinct geodesics
$G_0=G_0(\omega)$ and $G_1=G_1(\omega)$ such that
$$B_{G_0}(x,y) \neq B_{G_1}(x,y).$$
From this point a short argument allows us to conclude that
coexistence is possible with positive probability.

\section{Proof}
The heart of the proof is applying the ergodic theorem to the
\booze function. This is done in Lemmas \ref{lem4} and \ref{lem1}.
We start by showing that the symmetry of $\fpp$ implies that the
expected value of the \booze function is 0.
\begin{lemma}\label{lem0}
If $\{B(x,y)\}_{x,y \in \zd}$ is well defined a.s.\ then for all
$v \in \zd$
$$\E(B(\0,v))=0.$$
\end{lemma}

\begin{proof}
By symmetry of $\fpp$ we have that $\E(B(\0,\ez))=\E(B(\ez,\0))$.
Combining this with the fact that $B(\0,\ez)+B(\ez,\0)=0$ proves
the lemma.
\end{proof}

Now we apply the ergodic theorem to $B(\0,v)$.
\begin{lemma}\label{lem4}
If $\{B(x,y)\}_{x,y \in \zd}$ is well defined a.s.\ then for all
$v \in \zd$ and $\epsilon>0$ there exists $M$ such that
$$P\left(|B(\0,mv)|<\epsilon m \mbox{ for all $m>M$} \right)>1-\epsilon.$$
\end{lemma}
\begin{proof}
First rewrite $B(\0,mv)$ as follows.
\begin{eqnarray}
B(\0,mv)& =& B(\0,v)+B(v,2v) + \dots +B((m-1)v,mv) \nonumber\\
B(\0,mv)& =& B^{\omega}(\0,v)+ B^{T^v(\omega)}(\0,v)+ \dots +B^{T^{(m-1)v}(\omega)}(\0,v) \nonumber\\
B(\0,mv)& =& \sum_{j=0}^{m-1}B^{T^{jv}(\omega)}(\0,v)
\label{line1}
\end{eqnarray}
As $\fpp$ is good it is totally ergodic and the action
$(T^v,\fpp)$ is ergodic.  Thus by line (\ref{line1}) and Lemma
\ref{lem0} the claim is a consequence of the ergodic theorem.
\end{proof}

We now strengthen this lemma by using the following corollary of
shape theorem.  For $x \in \zd$ we let
$|x|=|x_1|+|x_2|+\dots+|x_d|.$

\begin{cor} \label{shape2}
There exist $0<k_1<k_2<\infty$ such that for every $\epsilon>0$
there exists an $N$ such that
$$P\left(k_1<\frac{\tau(\0,x)}{|x|}<k_2
 \mbox{ for all $x$ such that $|x|>N$}\right)>1-\epsilon.$$
\end{cor}
\begin{proof}
The existence of $k_2$ is due to the fact that the set $S$ (from
Theorem \ref{shape1}) is nonempty. The existence of $k_1$ follows
because one of the requirements of $\fpp$ being good is that $S$
is bounded.
\end{proof}

\begin{lemma} \label{lem1}
If $\{B(x,y)\}_{x,y \in \zd}$ is well defined a.s.\ then for any
$\epsilon>0$ there exists $N$ such that if $n>N$ then
$$P\left(\frac{B(\0,x)}{|x|}<\epsilon
 \mbox{ for all $x$ such that  $|x|=n$}\right)>1-\epsilon.$$
\end{lemma}

\begin{proof}
Given $\epsilon >0$ pick vectors $v_1,v_2,\dots, v_j$ such that
$|v_1|=|v_2|=\dots=|v_j|$ and for all $x$ sufficiently large there
exists $i \in \{1,2,\dots,j\}$ and $m \in \N$ such that
$$|x-mv_i|<\epsilon  |x| \mbox{ and } m|v_i| \leq |x|.$$
For all $x$ and $y$ we have that
$$B(\0,x)=B(\0,y)+B(y,x).$$
This implies that for any $x$ and $y$
$$B(\0,x)\leq B(\0,y)+\tau(y,x).$$
For any $n$ let $m$ be the largest integer such that $m|v_i|\leq
n$.  (This is independent of $i$.) Thus if there exists $x$ with
$|x|=n$ and $\frac{B(\0,x)}{|x|} \geq \epsilon$ then there exists
$i$ such that either
\begin{enumerate}
\item $B(\0,mv_i) \geq \epsilon n/2=\epsilon |x|/2$, or
\item $|x-mv_i|<\epsilon |x|/2k_2$ and $\tau(x,mv_i) \geq \epsilon |x|/2$.
\end{enumerate}
(The constant $k_2$ is from Corollary \ref{shape2}.)

By Lemma \ref{lem4} there exists $M$ such that
$$P(\mbox{there exists $m>M$ and $i \in\{1,2,\dots,j\}$ such that
$B(\0,mv_i)>2\epsilon m |v_i|/3>\epsilon n/2$}) <2\epsilon/3.$$
Thus the probability of the first event is less than $2\epsilon/3$
if $n$ is sufficiently large.

By Corollary \ref{shape2} there exists $L$ such that for any $l>L$
$$P(\mbox{there exists $z$ with $|z| \leq l$ and $\tau(\0,z)\geq k_2l$})<\epsilon/3j .$$
Applying this with each $mv_i$ in place of $\0$ and $\epsilon
n/2k_2$ in place of $l$ we get that
 the probability of the second event is less than $\epsilon/3$ if $n$ is sufficiently large.
Thus for any $\epsilon>0$ we get $N$ so that  if $n>N$ we get that
$$P\left(\mbox{there exists $x$ such that $|x|=n$ and } \frac{B(\0,x)}{|x|} \geq \epsilon
\right) < \epsilon$$ which proves the lemma.
\end{proof}

Next we show that this generates a contradiction with the shape
theorem.

\begin{lemma}
\label{coalesce} 
$P(\mbox{$\{B(x,y)\}_{x,y \in \zd}$ is well defined})=0.$
\end{lemma}

\begin{proof} We work by contradiction. Suppose that with positive
probability $\{B(x,y)\}_{x,y \in \zd}$ is well defined.  The
\booze function being well defined is a shift invariant event
which, by the ergodicity of $\fpp$, implies that $\{B(x,y)\}_{x,y
\in \zd}$ is well defined a.s.\ and the conclusions of Lemma
\ref{lem1} apply. Pick $\epsilon<\frac13 \min(k_1,1)$, where $k_1$
comes from Corollary \ref{shape2}.  By the choice of $\epsilon$
and Corollary \ref{shape2} we have that there exists $N$ such that
for all $n>N$ \be P\left(\frac{\tau(\0,x)}{|x|}>2\epsilon
 \mbox{ for all $x$ such that $|x|=n$}\right)>\frac23. \label{line2} \ee
By Lemma \ref{lem1} there exists $n>N$ such that
 \be P\left(\frac{B(\0,x)}{|x|}<\epsilon
 \mbox{ for all $x$ such that $|x|=n$}\right)>\frac23. \label{line3}\ee
However there exists at least one infinite geodesic
$G=(\0,v_1,v_2,\dots)$ which begins at $\0$. (The choice of $G$ is
immaterial.) For all $n$ there exists $k$ such that $|v_k|=n$. For
any $k$ we have that $B(\0,v_k)=\tau(\0,v_k)$. This shows that
lines (\ref{line2}) and (\ref{line3}) cannot both be true. Thus
the lemma is proven.
\end{proof}

%
%


Note that the lack of a well defined \booze function implies that
there exists at least two disjoint infinite geodesics.  Now we
show that the lack of a well defined \booze function also implies
coexistence has positive probability. Coexistence is implied if
there exist two infinite geodesics $G_0=(v_0,v_1,v_2,\dots)$ and
$G_1=(w_0,w_1,w_2,\dots)$ such that
$$B_{G_0}(\0,\ez) < 0< B_{G_1}(\0,\ez).$$
We show coexistence is possible by showing that we have two such
geodesics with positive probability.

%

 \noindent
 {\bf Proof of Theorem \ref{main}:\ }
By Lemma \ref{coalesce} we get an event $\tilde A$ of positive
probability and $x,y \in \zd$ such that for all $\omega \in \tilde
A$ we have two geodesics $G_0=G_0(\omega)=(v_0,v_1,v_2,\dots)$ and
$G_1=G_1(\omega)=(w_0,w_1,w_2,\dots)$ with
$$B_{G_0}(x,y) < B_{G_1}(x,y).$$
(If there is more than one pair of geodesics which satisfy this
equation we can choose $G_0$ and $G_1$ in any measurable manner.)
It causes no loss of generality to assume that $|x-y|=1$.  Thus by
the symmetry of $\fpp$ we can assume that $x=\0$ and $y=\ez$. As
$B_{G_0}(\0,\ez)$ and $B_{G_1}(\0,\ez)$ do not depend on any
finite number of edges in the geodesics, it causes no loss of
generality to assume that $\0,\ez$ are not endpoints of any of the
edges in  $G_0$ or $G_1$. By restricting to a smaller event
$A\subset \tilde A$ of positive probability we get a nonrandom
$r>0$ such that for all $\omega\in A$
\begin{equation} B_{G_0}(\0,\ez)<  r <B_{G_1}(\0,\ez). \label{g0} \end{equation}
By the symmetry of $\fpp$ we can assume $r\geq 0$. From the
definition of $B_{G_1}(\0,\ez)$ we get that $B_{G_1}(\0,\ez) \leq
\tau(\0,\ez)$.


Now we form a new event $A'$. Given $\omega \in A$ define
$\omega'$ by
$$\omega'(v,w)=\left\{%
\begin{array}{ll}
    \omega(v,w)+r, & \hbox{if $\ez \in \{v,w\}$;} \\
    \omega(v,w), & \hbox{else.}
\end{array}%
\right.$$ The event $A'$ consists of all  $\omega'$ that can be
formed from in this way from some $\omega \in A$. By conditions
\ref{cont0}, \ref{cont1} and \ref{cont2} of the definition of
$\mu$ being good, the event $A'$ also has positive measure.  We
let $\tau'$ indicate the passage times in $\omega'$ and $\tau$
indicate the passage times in $\omega$. It is easy to check that
for any $z \neq \ez$
$$\tau'(\ez,z)=\tau(\ez,z)+r.$$  Also if $\ez$ is not an endpoint of any of the edges in
the geodesic $G^{\omega}(\0,z)$ then
$$\tau'(\0,z)=\tau(\0,z).$$
As $B_{G_0}(\0,\ez)<B_{G_1}(\0,\ez)\leq \tau(\0,\ez)$ we have that
for all large $n$ the vertex $\ez$ is not an endpoint of any of
the edges in the geodesic $G^{\omega}(\0,v_n)$. Thus
$\tau'(\0,v_n)=\tau(\0,v_n)$ for all large $n$. Also note that
since neither $\0$ or $\ez$ is an endpoint of any of the edges
$G_0(\omega)$ or $G_1(\omega)$ we have that $G_0(\omega)$ and
$G_1(\omega)$ are both geodesics for $\omega'$.

Thus for any $\omega' \in A'$ we have that
\begin{eqnarray*}
 B^{\omega'}_{G_0(\omega)}(\0,\ez)&=&\lim_{n \to \infty} \tau'(\0,v_n)-\tau'(\ez,v_n)\\
 &=&\lim_{n \to \infty} \tau(\0,v_n)-(\tau(\ez,v_n)+r)\\
 &=&B^{\omega}_{G_0(\omega)}(\0,\ez)-r\\
 &<&0.
\end{eqnarray*}
The last step follows from line (\ref{g0}). We also get that
\begin{eqnarray*}
 B^{\omega'}_{G_1(\omega)}(\0,\ez)&=&\lim_{n \to \infty} \tau'(\0,w_n)-\tau'(\ez,w_n)\\
 &\geq&\lim_{n \to \infty} \tau(\0,w_n)-(\tau(\ez,w_n)+r)\\
 &\geq&B^{\omega}_{G_1(\omega)}(\0,\ez)-r\\
 &>&0.
\end{eqnarray*}
The last step follows from line (\ref{g0}). Thus we have
coexistence for all $\omega' \in A'$.

\section*{Acknowledgments}
The author thanks Itai Benjamini, Yuval Peres, and Oded Schramm
for helpful conversations.  He also thanks an anonymous referee
for pointing out an error in a previous version.

 \qed

\end{document}